\documentclass[12pt]{article}

\usepackage{amsmath}
\usepackage{amsthm}
\usepackage{amsfonts}

\begin{document}

\title{Elementary fixed points of the BRW smoothing transforms with infinite number of
summands}\date{}
\author{Aleksander M. Iksanov\footnote {e-mail address:
iksan@unicyb.kiev.ua} \\ \small{\emph{Faculty of Cybernetics},
\emph{Kiev Taras Shevchenko National University}},\\
\small{\emph{01033 Kiev, Ukraine}}}\maketitle
\begin{abstract}
The branching random walk (BRW) smoothing transform $T$ is defined
as $T:\text{distr}(U_{1})\mapsto \text{distr} \left(
\sum_{i=1}^{L}X_{i}U_{i}\right)$, where given realizations
$\{X_{i}\}_{i=1}^{L}$ of a point process, $U_{1},U_{2},\ldots$ are
conditionally independent identically distributed random
variables, and $0\leq \text{Prob}\{L=\infty \}\leq 1$. Given
$\alpha \in (0,1]$, $\alpha$-\emph{elementary} fixed points are
fixed points of $T$ whose Laplace-Stieltjes transforms $\varphi$
satisfy $\underset{s\rightarrow
+0}{\lim}\dfrac{1-\varphi(s)}{s^{\alpha}}=m$, where $m$ is any
given positive number. If $\alpha=1$, these are the fixed points
with finite mean. We show exactly when elementary fixed points
exist. In this case these are the only fixed points of $T$ and are
unique up to a multiplicative constant. These results do not need
any moment conditions. In particular, Biggins' martingale
convergence theorem is proved in full generality. Essentially we
apply recent results due to Lyons (1997) and Goldie and Maller
(2000) as the key point of our approach is a close connection
between fixed points with finite mean and perpetuities. As a
by-product, we lift from our general results the solution to a
Pitman-Yor problem. Finally, we study the tail behaviour of some
fixed points with finite mean.


\end{abstract}
\newpage \thispagestyle{plain}

\section{Introduction}

Unless otherwise stated, all random variables (rvs) studied in the
paper are assumed to be defined on a fixed probability space
$(\Omega ,\mathcal{F},\mathbb{P})$. We also assume that this
probability space is large enough to accomodate independent copies
of rvs. All distributions mentioned below will be probability
ones. Therefore the adjective will be usually dropped. The
distribution of an rv $X=X(\omega ) $, $\omega \in \Omega $ will
be denoted by $\mathcal{L}(X)$, and the degenerate distribution at
$x\geq 0$ (the delta measure) will be denoted by $\delta _{x}$.
Furthermore, $\mathcal{P}^{+} $ denotes the set of all Borel
probability measures on the nonnegative half line
$\mathbb{R}^{+}=[0,\infty )$.

Let $Z(\cdot)$ be a point process on $[0,+\infty )$, i.e. a
random, locally finite on $(0,+\infty)$, counting measure. It is
assumed that realizations (points) $\{X_{i}\}_{i=1}^{L}$ of $Z$
constitute a nonincreasing
collection of $L$ nonzero rvs, where $L$ is an rv with $%
\mathbb{P}\{L=\infty \}\in \lbrack 0,1]$. We consider the ordered
collection just for notational convenience, and this does not
restrict generality. Note that $\mathbb{E}Z[0,t]$ and
$\mathbb{E}Z(t,+\infty)$ may be infinite. Hence, the intensity
measure $\chi $ of $Z$ is defined for some positive finite $A$ as
follows $\int_{A}^{t+}\chi(dz)=\mathbb{E}Z(A,t]$, if $t\geq A$,
and $\int_{t-}^{A}\chi(dz)=\mathbb{E}Z(t,A]$, otherwise. Note that
in general  $\chi $ is a $\sigma $-finite Borel measure.

Let us now recall what\emph{\ the branching random walk (BRW)} is.
Assume that an initial ancestor is placed at the origin of the
real line and after one unit of time she gives birth to children
who form the first generation. Their displacements from the origin
are given by the point process $Z^{(1)}(B):=Z(e^{-B})$, where $B$
is a Borel set and $e^{-B}=\{e^{-x}:x\in B\}$, with points
$\{-\log X_{i}\}_{i=1}^{L}$. Each of these children also lives one
unit of time and has offspring in a like manner, so that the
positions of each family relative to the parent are given by an
independent copy of the point process $Z^{(1)}$. All children born
to individuals of the first generation forms the second generation
with
positions given by the point process $Z^{(2)}$ and so on. Thus $%
Z^{(n)}$ is the $n$-th generation point process. The discrete time
process $Z^{(0)}(B):=1_{\{\{0\}\in B\}}$ a.s., $Z^{(n)}$,
$n=1,2,\ldots$, is called the BRW. \newline Let $\mathcal{F}_{n}$
be the $\sigma $-fields containing all information about the first
$n$ generations, $n=1,2,\ldots$. It is well-known that, when the
mean number $m$ of children born to a person satisfies $m\in
(1,\infty]$, and $m(\gamma ):=\mathbb{E}\int_{-\infty }^{\infty
}e^{-\gamma t}Z^{(1)}(dt)\in (0,\infty )$ , for some $\gamma\geq
0$,
\begin{equation}
W^{(n)}(\gamma )=(1/m^{n}(\gamma ))\int_{-\infty }^{\infty
}e^{-\gamma t}Z^{(n)}(dt)
\end{equation}
is a nonnegative martingale with respect to $\mathcal{F}_{n}$. For
more information on the BRW and associated martingales see, for
example, Biggins (1977), Biggins and Kyprianou (1997).

Let $t_{r}$ be a rooted family tree associated with a point
process $Z^{(1)}$. We say that $(t_{r},X)$ is a labelled tree if
each individual (vertex) $\theta\in t_{r}\backslash\{0\}$ is
assigned its displacement $X(\theta)$ from its parent. The BRW
defines a probability measure $\mu$ on the set of labelled trees.

We address the problem of the existence and uniqueness of \emph{\
special} distributions of the \emph{nonnegative} rvs $W$
satisfying the following distributional equality
\begin{equation}
W\overset{d}{=}\sum_{i=1}^{L}X_{i}W_{i}\text{,}
\end{equation}
where $W_{1},W_{2},\ldots$ are, conditionally on
$\{X_{i}\}_{i=1}^{L}$, independent copies of $W$. The equality (2)
is equivalent to
\begin{equation}
\varphi(s)=\mathbb{E}\prod_{i=1}^{L}\varphi (X_{i}s) \text {,}
\end{equation}
where $\varphi$ is the Laplace-Stieltjes transform (LST) of
$\mathcal{L}(W)$. \newline If $W$ satisfies (2), then it is
natural to refer to $\mathcal{L}(W)$ as the \emph{fixed point of
the (supercritical) branching random walk (the BRW) smoothing
transform}
\begin{equation*}
\ \mathbb{T}:\mathcal{P}^{+}\rightarrow \mathcal{P}^{+}\cup \{\delta_{\infty}\}\text{; }\mathcal{L}%
(U_{1})\mapsto \mathcal{L}\left( \sum_{i=1}^{L}X_{i}U_{i}\right)
\text{,}
\end{equation*}
where given $Z$, $U_{1},U_{2},\ldots$ are conditionally
independent identically distributed rvs.

The name is explained as follows. First, we have
$\mathbb{T}\mathcal{L}(W)=\mathcal{L}(W)$ (fixed point). Secondly,
the martingale $W^{(n)}(\gamma )$, with an appropriate $\gamma $,
either (a) converges in mean to the rv $W$ having unit mean, or
(b) $\underset{n\rightarrow \infty}{\lim} W^{(n)}(\gamma )=0$
almost surely; in this case, under some assumptions, there exists
a (Seneta-Heyde) norming $\{c_{n}\}$ which means that
$\underset{n\rightarrow \infty}{\lim} W^{(n)}(\gamma )/c_{n}=W$ in
distribution (properties of the BRW). The dichotomy (a-b)
regarding the limiting behaviour of $W^{(n)}(\gamma )$ is
justified by Lyons' (1997) change of measure construction (his
formula (2)) together with his formulae (5) and (6). The
Seneta-Heyde norming is not investigated here. We mention the
works Biggins and Kyprianou (1996, 1997) and Cohn (1997) where
this subject is studied for the supercritical BRW with
\begin{equation}
L<\infty \text{ almost surely (a.s.).}
\end{equation}
Of special interest, as indicated by the title of this paper, is
the case that
\begin{equation}
\mathbb{P}\{L=\infty \}>0\text{,}
\end{equation}
which, for the most part, we will concentrate on.  However, as
particular cases, we obtain previously known results proved under
assumption (4).

As should be clear from the title of the paper, we must introduce
a new notion of an \emph{elementary} and a \emph{nonelementary}
fixed point. Given $\alpha\in (0,1]$, we will say that a
distribution $\mu_{\alpha}$ is an $\alpha$-\emph{elementary fixed
point} of $\mathbb{T}$ if its LST $\varphi_{\alpha}$ satisfies
\begin{equation}
\underset{s\rightarrow
+0}{\lim}\dfrac{1-\varphi_{\alpha}(s)}{s^{\alpha}}=m\text{,}
\end{equation}
for some finite $m>0$. Note that a fixed point is $1$-elementary
if and only if it has finite mean. The set of \emph{elementary
fixed points} consists of all $\alpha$-elementary fixed points,
$\alpha\in (0,1]$. A fixed point will be called
\emph{nonelementary} if there is no $\alpha\in (0,1]$ for which it
is $\alpha$-elementary. Below we provide a rather full description
of the elementary fixed points. As the analysis of nonelementary
fixed points uses quite different arguments, results in that
direction will appear elsewhere.

The paper is organized as follows. In Section 2 we formulate our
main results with proofs deferred to Section 3. In Section 4 we
lift from our general results the solution to a so-called
Pitman-Yor problem. The Pitman-Yor problem is closely related to
the problem of the existence of fixed points of the so-called shot
noise transforms. These are obtained by putting in (2)
$X_{i}:=h(\tau _{i})$,where $h$ is a nonnegative Borel measurable
function and $\tau _{i}$, $i=1,2,\ldots$ is a Poisson flow. Some
comments and references are given in Section 5 and the paper
closes with an Appendix where some needed technical results are
collected.

In addition to the notation introduced above, other frequently
used notation includes:\newline (a) $\overline{\mu }$ is the
\emph{size-biased distribution} corresponding to a given
distribution $\mu $ with a finite mean; it is defined by the
equality
\begin{equation*}
\overline{\mu }(dx):=\left( \int_{0}^{\infty }y\mu (dy)\right)
^{-1}x\mu (dx);
\end{equation*}
if $Z$ is an rv with $\mathcal{L}(Z)=\mu $ then $\overline{Z}$ is
an rv with $\mathcal{L}(\overline{Z})=\overline{\mu }$;\newline
(b) given a $\sigma $-finite measure $M$ and $\gamma>0$, the
measure $M^{\ast }_{\gamma}$ is defined by $M^{\ast
}_{\gamma}(dx):=x^{\gamma}M(dx)$, it is convenient to put $M^{\ast
}:=M^{\ast }_{1}$;\newline (c) by a \emph{perpetuity} is meant an
rv $B_{1}+\sum_{i=2}^{\infty} A_{1}\ldots A_{i-1}B_{i}$, where
$(A_{i},B_{i})$, $i=1,2,\ldots$ are independent copies of a random
pair $(A,B)$;\newline (d) given the BRW smoothing transform
$\mathbb{T}$ and $\gamma\in (0,1)$, the \emph{modified transform}
$\mathbb{T}_{\gamma}$ is defined in the same way as $\mathbb{T}$
with the only difference being that the underlying point process
has points $\{X_{i}^{\gamma}\}_{i=1}^{L}$. Thus
\begin{equation*}
\mathbb{T_{\gamma}}:\mathcal{P}^{+}\rightarrow \mathcal{P}^{+}\cup \{\delta_{\infty}\}\text{; }\mathcal{L}%
(U_{1})\mapsto \mathcal{L}\left(
\sum_{i=1}^{L}X_{i}^{\gamma}U_{i}\right) \text{.}
\end{equation*}

For the reader's convenience, we would like to point out two
conventions to be in force throughout the paper.\newline
(\textbf{\ C1 }) Clearly, $\delta _{0}$ always satisfies (2).
Hence in what follows we will seek for other fixed-point
distributions, not indicating this explicitly.\newline (\textbf{\
C2 }) We will assume that the intensity measure $\chi$ satisfies
the equality $\chi ^{\ast }\{0\}=0$.

Requiring (C1) is only a matter of convenience. (C2) is only
needed in the proof of Proposition 1.

\section{Results}\subsection{The existence and
uniqueness} Our first statement gives necessary conditions for the
existence of \emph{arbitrary} fixed points of $\mathbb {T}$ and
asserts the regular variation of $1-\varphi(s)$ at zero, where
$\varphi$ is the LST of a fixed point. These results deal with
both cases (4) and (5) and were \emph{partially} known under (4)
and some additional moment restrictions. See Liu (1998, Theorems
1.1 and 1.2) for details. Note, however, that the validity of
Proposition 1 does not require a priori assumptions and hence the
result is relatively new even if (4) holds.\newline
\textbf{Proposition 1.} If there exists a fixed point with the LST
$\varphi $, then \newline a) there exist at most two values $\beta
_{1}\leq $ $\beta _{2}$, $\beta _{1} $, $\beta _{2}\in (0,1]$ and
at least one of these (take $\beta _{1}=\beta _{2}$ in the next
equality) such that \textbf{Condition D$_{\beta_{1}}$}, defined in
the line below,
\begin{equation*}
\mathbb{E}\sum_{i=1}^{L}X_{i}^{\beta _{k}}=1,\text{ }k=1,2
\end{equation*}
holds and
\begin{equation*}
\underset{n\rightarrow \infty }{\lim \inf }S_{n}^{(\beta
_{1})}=-\infty \text{ almost surely,}
\end{equation*}
where $S_{n}^{(\beta _{1})}$, $n=0,1,\ldots$ is the random walk:
\begin{equation*}
S_{0}^{(\beta _{1})}:=0 \text{, \ } S_{n}^{(\beta
_{1})}:=\sum_{j=1}^{n}\log B_{j}^{(\beta _{1})} \text{, \ }
n=1,2\ldots
\end{equation*}
and $B_{1}^{(\beta _{1})},B_{2}^{(\beta _{1})}\ldots$ are
independent copies of an rv $B^{(\beta _{1})}$ with\newline
$\mathcal{L}(B^{(\beta _{1})})=\chi ^{\ast }_{\beta_{1}}$;
\newline
b) $\underset{s\rightarrow +0}{\lim }\dfrac{1-\varphi (sz)}{1-\varphi (s)}%
=z^{\beta _{1}}$, $z\geq 0$.

In what follows we are considering elementary fixed points. From
the above Proposition, we know that the random walk $S_{n}^{(\beta
_{1})}$, $n=0,1,\ldots$ is oscillating or drifting to $-\infty$
(non-oscillating). One of the results of the next Proposition
below is that the elementary fixed points correspond to the
non-oscillating random walks $S_{n}^{(\beta _{1})}$,
$n=0,1,\ldots$ On the other hand, one may conjecture that
nonelementary fixed points could correspond to the random walks of
both types.

Theorem 2 contains the necessary and sufficient conditions for the
existence of elementary fixed points. Lyons (1997)
proved this assertion for fixed points with finite mean
(which are $1$-elementary fixed points in our terminology) under the side condition that $\mathbb{E}%
\sum_{i=1}^{L}X_{i}\log X_{i}$ is finite.

When Condition D$_{\beta_{1}}$ holds, notice that
$\mathcal{L}(B^{(\beta_{1})})=\chi ^{\ast }_{\beta_{1}}$ and set
$R_{\beta_{1}}:=\log B^{(\beta_{1})}$. For a distribution $\sigma
$, set
\begin{equation*}
I_{R_{\beta_{1}}}(\sigma ):=\int_{(1,\infty )}\dfrac{\log
x}{\int_{0}^{\log x}\mathbb{P}\{R_{\beta_{1}}\leq -y\}dy}\sigma
(dx)\text{.}
\end{equation*}
\textbf{Theorem 2}. For $\alpha\in (0,1]$, an $\alpha$-elementary
fixed point exists if and only if $\beta_{1}=\alpha$ and Condition
D$_{\beta_{1}}$ together with one of the next three conditions
holds
\newline{}(a) $-\infty <\mathbb{E}R_{\beta_{1}}<0$ and
$\mathbb{E}(\sum_{i=1}^{L}X_{i}^{\beta_{1}})\log
^{+}(\sum_{i=1}^{L}X_{i}^{\beta_{1}})<\infty $;\newline{}(b)
$\mathbb{E}R_{\beta_{1}}=-\infty $ and
$I_{R_{\beta_{1}}}(\overline{\mathcal{L}(\sum_{i=1}^{L}X_{i}^{\beta_{1}})})<\infty
$;\newline{}(c)
$\mathbb{E}R^{+}_{\beta_{1}}=\mathbb{E}R^{-}_{\beta_{1}}=+\infty
$, $I_{R_{\beta_{1}}}(\chi ^{\ast }_{\beta_{1}})<\infty $ and
$I_{R_{\beta_{1}}}(\overline{\mathcal{L}(\sum_{i=1}^{L}X_{i}^{\beta_{1}})})<\infty
$.\newline The conditions (a), (b), (c) are equivalent to the two
requirements:\newline (d) $\underset{ n \rightarrow \infty}{\lim}
S_{n}^{(\beta_{1})}=-\infty $ \text{ \ almost surely};\newline(e)
$I_{R_{\beta_{1}}}(\overline{\mathcal{L}(\sum_{i=1}^{L}X_{i}^{\beta_{1}})})<\infty
$.

Now we would like to reveal an idea of the proof.\newline (1)
\emph{Case $\alpha=1$}. Fixed points of $\mathbb{T}$ are scale
invariant. Thus it suffices to study fixed points with unit mean.
By Lemma 14, a fixed
point with unit mean exists if and only if the nonnegative martingale $W^{(n)}(\gamma)$, $%
n=1,2,\ldots $ given by (1) converges in mean to it. Therefore, it
follows from Lyons' (1997) change of measure construction that
such fixed points are closely connected with \emph{perpetuities}.
Once this relation has been realized, to deal with the existence
of these fixed points, we can use results on perpetuities from the
recent comprehensive treatment of Goldie and Maller (2000). Just
in this way, either of the conditions (a)-(c) of Theorem 2 ensures
the martingale convergence. Thus the case $\alpha=1$ of Theorem 2
can be viewed as a generalization of Biggins' (1977) martingale
convergence theorem.\newline (2) \emph{Case $\alpha\in (0,1)$}. As
soon as some results are available for $1$-elementary fixed
points, the corresponding statements for $\alpha$-elementary fixed
points are easily derived via the stable transformation. See the
proof of Theorem 2 (case $\alpha\in (0,1)$) for the precise
statement. Section 5 contains some references to other works
dealing with the stable transformation.

In the next Proposition we describe the set $\mathcal{H}$ of all
elementary fixed points. In fact, this set consists of the fixed
points with finite mean ($\alpha=1$) and the fixed points
($\alpha\in (0,1)$) obtained from the fixed points with finite
mean for the modified transform $\mathbb{T}_{\alpha}$ via the
stable transformation (7). As far as the uniqueness is concerned,
we show that provided $\mathcal{H}$ is nonempty, it coincides with
the set of all fixed points.\newline \textbf{Proposition 3}. Let
Condition D$_{\beta_{1}}$, (d) and (e) of Theorem 2 be valid. Set
$\alpha:=\beta_{1}$.
\newline{}(a) If $\alpha =1$ then, for each $m>0$, $\mathbb{T}$ has a
unique ($1$-elementary) fixed point with mean $m$. \newline (b) If
$0<\alpha <1$ then, for each $m>0$ in (6), $\mathbb{T}$ has a
unique $\alpha$-elementary fixed point $\mu_{\alpha}$ given by
\begin{equation}
\mu _{\alpha}(x,\infty )=\int_{0}^{\infty }s_{\alpha
}(xt^{-1/\alpha },\infty )\mu _{1}(dx),\ \ x>0,
\end{equation}
where $s_{\alpha }$ is the strictly stable positive distribution
with the index of stability $\alpha $, and $\mu _{1}$ is the fixed
point with mean $m$ of the modified transform
$\mathbb{T}_{\alpha}$.\newline (c) $\mathbb{T}$ on
$\mathcal{P}^{+} $ has no other fixed points than those described
in (a) and (b). \subsection{Tail behaviour} In this Section we
study the tail behaviour of fixed points with finite mean. First
we investigate the existence of moments of order $p>1$. As a
by-product we obtain conditions for the $L_{p}$-convergence of the
martingale $W^{(n)}(1)$. In case (4) the next Proposition is due
to Liu (2000, Theorem 2.3). Although his proof works well for the
infinite case too, we give an
alternative proof for the ''$\Rightarrow $'' part of the assertion and, when $%
p\in (1,2]$, for the ''$\Leftarrow $'' part.  \newline
\textbf{Proposition 4.} Assume that there exists a fixed point
$\mathcal{L}(W)\neq
\delta _{a}$, $a\geq 0$ with finite mean. Then, for each fixed $p>1$, $%
\mathbb{E}W^{p}<\infty $ if and only if
\begin{equation}
\mathbb{E}\left( \sum_{i=1}^{L}X_{i}\right) ^{p}<\infty \text{ and }\mathbb{E%
}\sum_{i=1}^{L}X_{i}^{p}<1.\newline
\end{equation}

The next result is obvious. In fact, it suffices to note that if $\mathbb{E}%
W^{p}<\infty $ then $W^{(n)}(1)=\mathbb{E}(\underset{m\rightarrow \infty }{%
\lim \inf }W^{(m)}(1)|\mathcal{F}_{n})$ and use Jensen's
inequality to see that $W^{(n)}(1)$ is bounded in $L_{p}$.\newline
\textbf{Corollary 5. }For each fixed $p>1$, the martingale
$W^{(n)}(1)$ is $L_{p}$-convergent if and only if (8) and
Condition D$_{1}$ hold.

Proposition 4 will be proved by using the modern technique based
on an appropriate change of measure. However, a simpler proof of
the ''$\Leftarrow $'' part is available as soon as one realizes
that under (8) the BRW smoothing transform $\mathbb{T}$ is a
strict contraction on some metric space. The reader may want to
consult R\"{o}sler (1992) and Rachev and R\"{u}schendorf (1995)
where the Contraction Principle is used to study transforms more
general than ours.

For fixed $\delta >1$ and $m>0$, let us consider the set $\mathcal{P}%
^{+}(\delta ,m)$ of distributions defined as follows
\begin{equation*}
\mathcal{P}^{+}(\delta ,m):=\{\mu \in
\mathcal{P}^{+}:\int_{0}^{\infty }x\mu (dx)=m,\int_{0}^{\infty
}x^{\delta }\mu (dx)<\infty \}\text{.}
\end{equation*}
Given $\mathcal{L}(Y)\in \mathcal{P}^{+}(\delta ,m)$ and a point
process whose points $\{X_{i}\}_{i=1}^{L}$ satisfy $\mathbb{E}%
(\sum_{i=1}^{L}X_{i})=1$ and $\mathbb{E}\left(
\sum_{i=1}^{L}X_{i}\right) ^{\delta }<\infty $, we have
$\mathcal{L}(\sum_{i=1}^{L}X_{i}Y_{i})\in
\mathcal{P}^{+}(\delta ,m)$, where $Y_{1},Y_{2},\ldots$ are, conditionally on $%
\{X_{i}\}_{i=1}^{L}$, independent copies of $Y$. Indeed, it is
easily seen that
\begin{equation*}
\mathbb{E}(\sum_{i=1}^{L}X_{i}Y_{i})=\mathbb{E}(\sum_{i=1}^{L}X_{i})\mathbb{E%
}Y=m\text{.}
\end{equation*}
Also by the convexity of the function $x\rightarrow x^{\delta }$,
we have
\begin{equation*}
\mathbb{E}(\sum_{i=1}^{L}X_{i}Y_{i})^{\delta }=\mathbb{E}(\mathbb{E}%
(\sum_{i=1}^{L}X_{i}Y_{i})^{\delta }/\mathcal{F}_{1}))\leq
\end{equation*}
\begin{equation*}
\leq \mathbb{E}(\mathbb{E(}(\sum_{i=1}^{L}X_{i})^{\delta
-1}(\sum_{i=1}^{L}X_{i}Y_{i}^{\delta }))/\mathcal{F}_{1}))=
\end{equation*}
\begin{equation*}
=\mathbb{E}(\sum_{i=1}^{L}X_{i})^{\delta
}\mathbb{E}Y^{\delta}<\infty .
\end{equation*}
Thus $\mathbb{T}$ maps $\mathcal{P}^{+}(\delta ,m)$ into itself.
\newline For $\mu _{1},\mu _{2}\in \mathcal{P}^{+}(\delta
,m)$, let us define the function
\begin{equation*}
r_{\delta }(\mu _{1},\mu _{2}):=\int_{0}^{\infty }s^{-\delta
-1}\left| \int_{0}^{\infty }\exp (isx)\mu
_{1}(dx)-\int_{0}^{\infty }\exp (isx)\mu _{2}(dx)\right|
ds\text{.}
\end{equation*}
From Lemma 3.1 of Baringhaus and Gr\"{u}bel (1997) we know that provided $%
\delta \in (1,2)$ $r_{\delta }$ is a metric on
$\mathcal{P}^{+}(\delta ,m)$ and that $(\mathcal{P}^{+}(\delta
,m),r_{\delta })$ is a complete metric space. \medskip \newline
\textbf{Proposition 6.} Let (8) with $p\in (1,2)$ and Condition
D$_{1}$ be in force. Then the BRW smoothing transform
$\mathbb{T}$, on $(\mathcal{P}^{+}(p,m),r_{p})$, is a strict
contraction. In particular, $\mathbb{E} W^{p}<\infty $.

The next assertion refines Theorem 2.2 of Liu (2000) and extends
its result to the transforms satisfying (5). First, we show that
the power-like tail behaviour does not depend on the type of
$\mathcal{L}(\log B^{(1)}).$ This features the fixed points under
consideration among general perpetuities (see Grincevi\v{c}ius
(1975), Theorem 2). Secondly, we indicate the explicit form of the
constant in the limit relation (9). \newline
\textbf{Proposition 7.} Assume that for some $b>1$, $\mathbb{E}%
(\sum_{i=1}^{L}X_{i})=\mathbb{E}(\sum_{i=1}^{L}X_{i}^{b})=1$, $\mathbb{%
E}(\sum_{i=1}^{L}X_{i}^{b}\log ^{+}X_{i})<\infty $ and $\mathbb{E}%
(\sum_{i=1}^{L}X_{i})^{b}<\infty $. Then there exist a fixed point
$\mu=\mathcal{L}(W)$ having finite mean, and a positive constant
$C_{b}$ such that
\begin{equation}
\underset{x\rightarrow \infty }{\lim }x^{b}\mu (x,\infty )=C_{b}
\text{.}
\end{equation}
Furthermore, \newline 1) if $\mathcal{L}(\log B^{(1)})$ is
nonarithmetic then
\begin{equation*}
C_{b}=(\mathbb{E}(\sum_{i=1}^{L}X_{i}^{b}\log
X_{i}))^{-1}\int_{0}^{\infty }y^{b-1}(\mu (y,\infty )-N(y,\infty ))dy%
\text{;}
\end{equation*}
\newline
where $N$ is the $\sigma $-finite measure defined by $N^{\ast }:=\mathcal{L}%
(B^{(1)}\overline{W})$. \newline 2) if $\mathcal{L}(\log B^{(1)})$
is arithmetic with the span $\varsigma $ then
\begin{equation*}
C_{b}=(\mathbb{E}(\sum_{i=1}^{L}X_{i}^{b}\log
X_{i}))^{-1}\sum_{k=-\infty }^{\infty }e^{-\varsigma
k}\int_{0}^{\exp \varsigma k}y^{b}(\mu (y,\infty )-N(y))dy\text{.}
\end{equation*}

\section{Proofs of the results} \textbf{Proof of
Proposition 1.} We consider the case (5), as the other one, when
(4) holds, can be treated similarly. Certainly, we will comment on
all points which require different arguments for these cases.
\newline Put $\psi (s):=\dfrac{1-\varphi (s)}{s}.$ From (3) we deduce that
\begin{equation}
1=\underset{s\rightarrow +0}{\lim }\mathbb{E}\sum_{i=1}^{L }X_{i}\dfrac{%
\psi (X_{i}s)}{\psi (s)}\prod_{k=1}^{i-1}\varphi (X_{k}s)\text{,}
\end{equation}
the empty product is always taken to be equal to $1.$ 
Further we will use arguments given in the proof
of Lemma 3.3 in Iksanov and Jurek (2002).\newline Since
\begin{equation}
0<\dfrac{\psi (sz)}{\psi (s)}\leq 1\ \ \text{for all}\ \ z\geq 1,
\end{equation}
by the selection principle, for any positive sequence $s_{n}$
which tends to $0$ as $n\rightarrow \infty $, there exists a
subsequence $s_{m_{n}}$ such that, for $t_{n}:=s_{m_{n}}$ and
$z>1$, $\dfrac{\psi (t_{n}z)}{\psi (t_{n})}$ converges to some
finite limit $\Lambda (z)$ as $n\rightarrow \infty $. On the other
hand, since each $\psi (t_{n}z)$ for $n=1,2,\ldots$, is a
completely monotone function in $z\in (0,\infty )$, and this
property is preserved under the limits, $\Lambda (z)$ is also
completely monotone, and thus, in particular, it is continuous on
$(0,\infty )$. Furthermore,
\begin{equation}
\underset{n\rightarrow \infty }{\lim }\dfrac{\psi (t_{n}z)}{\psi (t_{n})}%
=\Lambda (z)\text{ locally uniformly on }(0,\infty )\text{.}
\end{equation}
Also, for fixed $v>0$, we have that
\begin{equation}
\underset{n\rightarrow \infty }{\lim }\dfrac{\psi (t_{n}vz)}{\psi (t_{n}v)}=%
\dfrac{\Lambda (vz)}{\Lambda (v)}\text{ \ locally uniformly in }z\in (0,\infty )%
\text{.}
\end{equation}
If we would know that $\Lambda (\infty )=0$, the convergence in
(13) was uniform outside $0$, and we merely interchanged the limit
and the expectation. Under the current circumstances, in view of
(10) and Fatou's lemma, we obtain
\begin{equation}
1=\underset{n\rightarrow \infty }{\lim }\mathbb{E}\sum_{i=1}^{L }X_{i}%
\dfrac{\psi (X_{i}t_{n}v)}{\psi (t_{n}v)}\prod_{k=1}^{i-1}\varphi
(X_{k}s)\geq \mathbb{E}\sum_{i=1}^{L }X_{i}\dfrac{\Lambda (X_{i}v)}{%
\Lambda (v)}=
\end{equation}
\begin{equation*}
=\int_{0}^{\infty }\dfrac{\Lambda (vz)}{\Lambda (v)}\chi ^{\ast
}(dz)=:q\in (0,1]\text{, }v>0\newline \text{. }
\end{equation*}
\newline
After rewriting this in a more convenient form we get
\begin{equation*}
\int_{0}^{\infty }\Lambda (vz)(q^{-1}\chi ^{\ast }(dz))=\Lambda (v),v>0%
\newline
\text{.}
\end{equation*}
Changing of variable $z:=e^{-u}$ gives the integrated Cauchy
functional equation (in $\Lambda (e^{-\nu })$). It is known (see,
for example, Theorem 8.1.6 in Ramachandran and Lau (1991)) that
the solutions to such an equation are of the form
\begin{equation}
\Lambda (v)=p_{1}(v)v^{\beta _{1}-1}+p_{2}(v)v^{\beta _{2}-1},\
\text{\ for almost all }v>0\newline ;\
\end{equation}
\begin{equation}
\ p_{k}(v)=p_{k}(vw)\geq 0\ \ \text{for all}\ \ w\in supp(\chi
),k=1,2,
\end{equation}
where $\beta _{1}\leq \beta _{2}$ are determined by the equation
\begin{equation}
q=\int_{0}^{\infty }z^{\beta _{k}}\chi (dz)=\mathbb{E}\sum_{i=1}^{L}X_{i}^{%
\beta _{k}}\text{, }k=1,2.
\end{equation}
\newline
In our case, in view of continuity of $\Lambda $, (15) holds for
all $v>0$. [Note also that nonzero functions $p_{k}$ may be
different from the identical one only if
\begin{equation*}
supp(\chi )=\{\exp (n\gamma ):n\in \mathbb{Z}\},\text{ }
\end{equation*}
\begin{equation*}
\text{where }\gamma \text{ is the largest value that permits such
a representation.]\newline }
\end{equation*}
Now we must have $\beta _{1}\leq 1$, as otherwise $\Lambda $ given
by (15) is nondecreasing for some $\nu >0$. By the same reasoning,
if $\beta _{2}>1$ then $p_{2}(v)\equiv 0$. Note further that
$\beta _{2}$ cannot be negative, as $v\Lambda (v)$ being the
nonincreasing function for small enough $\nu $, would be the limit
of \emph{nondecreasing} functions $\dfrac{1-\varphi
(t_{n}v)}{1-\varphi (t_{n})}$. Finally, the case $\beta _{i}=0$,
$i=1,2$ is
excluded by (17). Indeed, if either (5) or both (4) and $%
\mathbb{E}L=\infty $ holds then the integral in (17) is infinite.
On the other hand, if $\mathbb{E}L<\infty $ then (17) implies that
$\mathbb{E}L\leq 1 $ which is impossible by Lemma 10(b) below. All
in all, it remains to consider two cases which we will study
separately: (A1) $0<\beta _{1}=\beta _{2}\leq 1$ (which simply
means that there is a unique $0<\beta \leq 1$ satisfying (17);
note however that there may be a $\beta >1$ satisfying (17)) and
(A2) $0<\beta _{1}<\beta _{2}\leq 1$. \

(A1) $\Lambda (v)=p(v)v^{\beta _{1}-1}$, where $p(v)\geq 0$
satisfies (16). We again may repeat the part of the proof of the
Lemma 3.3 in Iksanov and Jurek (2002) to conclude that $p(v)\equiv
1$ or equivalently $\Lambda (v)=v^{\beta _{1}-1}$. \newline It is
worth recording here as we need to use twice these arguments. Let
us introduce $k(v):=p(-\log v)$. Then $k(v)=v^{1-\beta
_{1}}\Lambda (v)$ is differentiable on $\mathbb{R}$ and \newline
\begin{equation*}
k^{^{\prime }}(v)=v^{1-\beta _{1}}((1-\beta _{1})v^{-1}\Lambda
(v)+\Lambda ^{^{\prime }}(v)).
\end{equation*}

Because of the differentiability and periodicity of $k(v)$ there exists a $%
v_{0}>0$ such that $k^{^{\prime }}(v_{0})=0$. In fact,
$k^{^{\prime }}(u^{n}v_{0})=0$, for $u\in supp(\chi )$ and
$n=1,2,\ldots$ On the other hand, both functions $v^{-1}\Lambda
(v)$ and $-\Lambda ^{^{\prime }}(v)$ are positive, nonincreasing
and convex. Consequently, for $0\leq \beta_{1} \leq 1$, the
equation
\newline
\begin{equation*}
(1-\beta _{1})v^{-1}\Lambda (v)=-\Lambda ^{^{\prime }}(v)
\end{equation*}
either holds identically or has at most two solutions (graphs of
the left- and the right-hand side may either coincide or intersect
at most at two points). However, the latter means that
$k^{^{\prime }}(v)=0$ at most at two
points, which contradicts the fact that $k^{^{\prime }}(u^{n}v_{0})=0$ for $%
n=1,2,\ldots$ Thus, $(1-\beta _{1})v^{-1}\Lambda (v)\equiv
-\Lambda ^{^{\prime
}}(v)$ which implies that $k(v)$ is a constant. Since by (12), $\Lambda (1)=1$%
, we conclude $k(v)=1$, for $v\geq 0$, or equivalently $\Lambda
(v)=v^{\beta _{1}-1}$.

With this $\Lambda $, the convergence in (13) is uniform outside $0$. In view of (C2) we have $%
\chi ^{\ast }\{0\}=0$. Therefore, we may interchange the limit and
the expectation in (10) to obtain:
\begin{equation*}
\mathbb{E}\sum_{i=1}^{L }X_{i}^{\beta _{1}}=1\text{.}
\end{equation*}
Furthermore, appealing to (12) we get
\begin{equation*}
\underset{n\rightarrow \infty }{\lim }\dfrac{1-\varphi
(t_{n}v)}{1-\varphi (t_{n})}=v\Lambda (v)=v^{\beta _{1}},\ \
\text{for all}\ \ v\geq 0.
\end{equation*}
However, the same argument below (11) can be repeated for any
subsequence, therefore we conclude that
\begin{equation*}
\underset{s\rightarrow +0}{\lim }\dfrac{1-\varphi (sz)}{1-\varphi (s)}%
=z^{\beta _{1}}\text{,}\ \ \text{for all}\ \ z\geq 0.
\end{equation*}
\newline
Thus in case (A1) the proof of b) and the first part of a) is
complete. To check the remaining part of a), set $\Phi _{\beta
_{1}}(s):=e^{\beta
_{1}s}(1-\varphi (e^{-s}))$ and $\Psi _{\beta _{1}}(s):=e^{\beta _{1}s}%
\mathbb{E}\{\prod_{i=1}^{L }\varphi
(e^{-s}X_{i})-1+\sum_{i=1}^{L}(1-\varphi (e^{-s}X_{i}))\}$. We
have the (infinite) analogue of the renewal equation given of
Lemma 2.3 by Durrett and Liggett (1983)
\begin{equation}
\Phi _{\beta _{1}}(t)=\mathbb{E}\Phi _{\beta _{1}}(t-S_{n}^{(\beta
_{1})})-\sum_{i=0}^{n-1}\mathbb{E}\Psi _{\beta
_{1}}(t-S_{i}^{(\beta _{1})}).
\end{equation}
In view of Lemma 10(a) (see Appendix), the random walk
$S_{n}^{(\beta _{1})},n=0,1,\ldots$ is
nondegenerate at $0$. So let us assume that $\underset{n\rightarrow \infty }{%
\lim \sup } S_{n}^{(\beta _{1})}=+\infty $ a.s. Thus there exists
a nonrandom subsequence $\{n_{k}\}$ which approaches
infinity together with $k$, and such that $\underset{k\rightarrow \infty }{%
\lim } S_{n_{k}}^{(\beta _{1})}=+\infty $ a.s. As $\Phi _{\beta
_{1}}$ is bounded on the neighborhoods of infinity and
$\underset{s\rightarrow -\infty }{\lim }\Phi
_{\beta _{1}}(s)=0$, we have $\underset{k\rightarrow \infty }{\lim }\mathbb{E}%
\Phi _{\beta _{1}}(S_{n_{k}}^{(\beta _{1})})=0$. Therefore,
$\underset{k\rightarrow \infty
}{\lim}\sum_{i=0}^{n_{k}}\mathbb{E}\Psi _{\beta
_{1}}(S_{i}^{(\beta _{1})})$ exists and it is nonpositive and
possibly infinite. This contradicts the fact that $\Phi _{\beta
_{1}}(0)=1-\varphi (1)>0$. [Note that Biggins (1977, Lemma 3) used
similar arguments]. This concludes the proof of the Proposition in
case (A1).

(A2) In this case there exists $1<w\in supp(\chi )$. From (15) and
(16) we have, for any $n\in \mathbb{N}$,
\begin{equation*}
\Lambda (w^{n})=p_{1}(w^{n})(w^{n})^{\beta
_{1}-1}+p_{2}(w^{n})(w^{n})^{\beta _{2}-1}=p_{1}(1)(w^{n})^{\beta
_{1}-1}+p_{2}(1)(w^{n})^{\beta _{2}-1}\text{.}
\end{equation*}
If $\beta _{2}<1$, this yields $\underset{\nu \rightarrow \infty
}{\lim }\Lambda (\nu )=0$, in view of the monotonicity.
Consequently, in the same way as in the study of case (A1) we have
in (17) $q=1$, thus proving that Condition D$_{\beta_{1}}$ holds.
If $\beta _{2}=1$, we just divide the sum in (14) into two parts:
$\lim \mathbb{E} \sum_{1}^{L}=\lim \mathbb{E}
\sum_{1}^{Z[0,d]}+\lim \mathbb{E} \sum_{Z[0,d]+1}^{L}$, for some
$d>0$. Now we may interchange the limit and either of sums
separately. While for the first sum we use the local uniformity of
convergence in (13), for the second one we use nonincreasingness
of $\psi $, the dominated convergence and the fact that $\chi
^{\ast}$ is the finite measure. This implies that $q=1$ and hence
Condition D$_{\beta_{1}}$ holds.

Now the distribution of the rv $B^{(\beta _{1})}$ is well-defined
and, moreover,\newline
$\mathbb{E}(B^{(\beta _{1})})^{\beta _{2}-\beta _{1}-\varepsilon }<1$, for $%
\varepsilon \in (0,\beta _{2}-\beta _{1})$. By making use of
Jensen's inequality, we conclude that $\mathbb{E}\log B^{(\beta
_{1})}<0$ which implies $\underset{n\rightarrow \infty }{\lim
}S_{n}^{(\beta _{1})}=-\infty $ a.s.\newline
We now turn to the proof of the regular variation. By (17), there exists a $%
1>w\in supp(\chi )$. Thus (15) and (16) give for any\ $n\in
\mathbb{N}$,
\begin{equation*}
(w^{n})^{1-\beta _{1}}\Lambda
(w^{n})=p_{1}(w^{n})+p_{2}(w^{n})(w^{n})^{\beta _{2}-\beta
_{1}}=p_{1}(1)+p_{2}(1)(w^{n})^{\beta _{2}-\beta _{1}}\text{.}
\end{equation*}
\newline
Consequently, we have $\underset{n\rightarrow \infty }{\lim
}(w^{n})^{1-\beta _{1}}\Lambda (w^{n})=p_{1}(1)$. Using the
monotonicity of $\Lambda $, we get
\begin{equation}
wp_{1}(1)\leq \underset{\nu \rightarrow +0}{\lim \inf } \nu
^{1-\beta _{1}}\Lambda (\nu )\leq \underset{\nu \rightarrow
+0}{\lim \sup } \nu ^{1-\beta _{1}}\Lambda (\nu )\leq
\frac{1}{w}p_{1}(1).
\end{equation}
\newline
In view of (18), $\Phi _{\beta _{1}}(t-S_{n}^{(\beta
_{1})})-\sum_{i=0}^{n-1}\Psi _{\beta _{1}}(t-S_{i}^{(\beta _{1})})$, $%
n=0,1,\ldots$ is a martingale. Note that here the functions $\Phi
_{\beta _{1}}$ and $\Psi _{\beta _{1}}$ are constructed in the
same way as above (18), but with $\varphi $ and $\{X_{i}\}$ that
we are currently studying. As the random walk $S_{n}^{(\beta
_{1})}$ drifts to $-\infty$, the stopping time $\tau =\min
\{n>0:S_{n}^{(\beta _{1})}<0\}$ is a.s. finite. By the martingale
stopping theorem, we have
\begin{equation*}
\Phi _{\beta _{1}}(t)=\mathbb{E}\Phi _{\beta _{1}}(t-S_{\tau
\wedge n}^{(\beta _{1})})-\sum_{i=0}^{\tau \wedge n-1}
\mathbb{E}\Psi _{\beta _{1}}(t-S_{i}^{(\beta _{1})})\leq
\mathbb{E}\Phi _{\beta _{1}}(t-S_{\tau\wedge n }^{(\beta _{1})}).
\end{equation*}
To get the latter inequality, we have used nonnegativity of $\Psi
_{\beta _{1}}$(see Durrett and Liggett (1983, Lemma 2.4b) for the
proof). Put $\Lambda _{\beta _{1}}(\nu ):=\nu ^{1-\beta
_{1}}\Lambda (\nu )$. Using the same $t_{n}$ as in (12) and the
local uniform convergence there gives
\begin{equation*}
1\leq \underset{m\rightarrow \infty }{\lim }\int_{0^{+}}^{1}z^{1-\beta _{1}}%
\dfrac{\psi (t_{m}vz)}{\psi (t_{m}v)}\mathcal{L}(e^{S_{\tau\wedge
n }^{(\beta _{1})}})(dz)=
\end{equation*}
\begin{equation*}
=\int_{0}^{1}\dfrac{\Lambda _{\beta _{1}}(\nu z)}{\Lambda _{\beta _{1}}(\nu )%
}\mathcal{L}(e^{S_{\tau\wedge n }^{(\beta _{1})}})(dz)<\infty.
\end{equation*}
The integrand is bounded according to (19). Hence the Lebesgue
bounded convergence allows us to pass to the limit as
$n\rightarrow \infty$, to get
\begin{equation*}
1\leq \int_{0}^{1}\dfrac{\Lambda _{\beta _{1}}(\nu z)}{\Lambda _{\beta _{1}}(\nu )%
}\mathcal{L}(e^{S_{\tau}^{(\beta _{1})}})(dz)=:Q<\infty.
\end{equation*}
Now we may repeat the discussion of the first part of the proof to
conclude
\begin{eqnarray*}
\Lambda _{\beta _{1}}(\nu ) &=&P(v)v^{b_{1}-1}\ \text{\ for all
}v>0\newline
; \\
\ \ P(v) &=&P(vw)\geq 0\ \ \text{for all}\ \ w\in supp(\mathcal{L}%
(e^{S_{\tau }^{(\beta _{1})}}))\text{.}
\end{eqnarray*}
where $b_{1}$, which is necessarily unique as $0<e^{S_{\tau
}^{(\beta _{1})}}<1$ a.s., is determined by the equation
\begin{equation*}
Q=\int_{0}^{1}z^{b_{1}-1}\mathcal{L}(e^{S_{\tau }^{(\beta
_{1})}})(dz).
\end{equation*}
\newline
Suppose that $b_{1}\neq 1$. Then $b_{1}<1$ which implies that
$\Lambda _{\beta _{1}}(\nu )$ is unbounded near zero. A
contradiction. Thus, $\Lambda (\nu )=P(v)v^{\beta _{1}-1}$, and it
remains to copy the proof of case (A1) to show that $P(v)\equiv
1$, $\nu \geq 0$. Repeating all arguments above for each
subsequence (like $t_{n})$ finishes the proof of the Proposition.
\newline \textbf{Proof of Theorem 2.} \emph{Case
$\alpha=1$}. We first notice that if $\mathcal{L}(W)$ satisfies
(2), so is $\mathcal{L}(cW)$, for any $c>0$. Thus it suffices to
study the situation when $\mathbb{E}W=1$. That the conditions
(a)-(c) and (d)-(e) are equivalent follows from Theorem 2.1 and
Corollary 4.1 of Goldie and Maller (2000).\newline The part (a)
was studied by Lyons (1997). As to (b), Lyons noticed that his
approach works well in this case too, but did not provide
arguments. Let us show simultaneously that each of the couples of
conditions (d), D$_{1}$ and (c), D$_{1}$ is sufficient.

Lyons (1997) constructed a probability space (($t,X,\xi $), $\mathcal{F}%
^{\ast }$, $\widehat{\mu }^{\ast }$)$,$ where ($t,X,\xi $) is a
space of
infinite labelled trees ($t,X$) with distinguished rays $\xi $, $\mathcal{F}%
^{\ast }=\cup \mathcal{F}_{n}^{\ast }$, where $\mathcal{F}_{n}^{\ast }$, $%
n=1,2,\ldots$ are the $\sigma $-fields containing all information
about the first $n$ generations in ($t,X,\xi $), and $\widehat{\mu
}^{\ast }$ is a probability measure whose ''double'' restriction
$\widehat{\mu }_{n}$, first to ($t,X$) then to $\mathcal{F}_{n}$
satisfies
\begin{equation*}
\frac{d\widehat{\mu }_{n}}{d\mu _{n}}=W^{(n)}(1)\text{, for all
}n\text{ and all }(t,X),
\end{equation*}
where $\mu _{n}$ is the restriction of $\mu $ to\
$\mathcal{F}_{n}$. Let $S$ be an rv whose distribution is given as
follows
\begin{equation*}
d\mathcal{L}(S)=(\sum_{i=1}^{L}X_{i})d\mathcal{L}(Z)\text{.}
\end{equation*}
Then with $\mathcal{G}$ being the $\sigma $-field generated by the
copies of
$S$ we have $V_{n}:=\mathbb{E}_{\widehat{\mu }^{\ast }}(W^{(n)}(1)/\mathcal{G%
})=:V_{1,n}-V_{2,n}\leq V_{1,n}$, where
$V_{1,n}=N_{1}+\sum_{k=1}^{n-1}M_{1} \ldots M_{k}N_{k+1}$ and
$V_{2,n}=\sum_{k=1}^{n-1}M_{1} \ldots M_{k}$. Here
$M_{1},M_{2},\ldots$ are $\widehat{\mu }^{\ast }$ iidrvs with the
distribution $\chi ^{\ast }$ which are also $\widehat{\mu }^{\ast
}$ independent of $\widehat{\mu }^{\ast
}$ iidrvs $N_{1},N_{2},\ldots$ with the distribution $\overline{\mathcal{L}%
(\sum_{i=1}^{L}X_{i})}$. Thus, in view of Fatou's lemma and Lemma
11, for the existence of $\mathcal{L}(W)$ it suffices that
$V_{1,n}$ be $\widehat{\mu }^{\ast }$ a.s. convergent (with a
limit being a perpetuity). While the condition (b) ensures the
convergence of $V_{1,n}$ by Corollary 4.1(c) of Goldie and Maller
(2000), (c) ensures this by Corollary 4.1(d) of the same
reference.

Let us now assume that $\mathcal{L}(W)$ exists. The necessity of
Condition D$_{1}$ follows from the equality $\mathbb{E}W=(\mathbb{E}W_{1})\mathbb{E}%
\sum_{i=1}^{L}X_{i}$. Therefore, in the sequel we may and do
assume that the
distributions $\chi ^{\ast }$ and $\overline{\mathcal{L}(\sum_{i=1}^{L}X_{i})%
}$ \ are well-defined. Using formula (7) of Lyons (1997)
\begin{equation}
W^{(n+1)}(1)\geq Y_{n+1}\text{, }n=0,1,\ldots\text{,}
\end{equation}
where $Y_{n+1}\overset{d}{=}M_{1} \ldots M_{n}N_{n+1}$ and our Lemma 5.2 with $%
\varepsilon =0$ (which is just the formula (6) of Lyons (1997))
gives
\begin{equation*}
\underset{k\rightarrow \infty }{\lim \sup }M_{1} \ldots M_{k}N_{k+1}<\infty ,%
\widehat{\mu }^{\ast }\text{a.s. }
\end{equation*}

Suppose that the condition $\underset{n\rightarrow \infty }{\lim}
S_{n}^{(1)}=-\infty$ $\widehat{\mu}^{\ast }$ a.s., fails to hold. This implies $%
\underset{k\rightarrow \infty }{\lim \sup }M_{1} \ldots M_{k}=\underset{%
k\rightarrow \infty }{\lim \sup }\exp S_{k}^{(1)}=+\infty $ $%
\widehat{\mu }^{\ast }$ a.s. A contradiction. Thus, in what
follows we assume that $\underset{k\rightarrow \infty }{\lim
}M_{1} \ldots M_{k}=0$ $\widehat{\mu} ^{\ast}$ a.s.

As $M_{1}1_{\{M_{1}\in \lbrack 0,1]\}}\ldots M_{k}1_{\{M_{k}\in
\lbrack
0,1]\}}\leq M_{1}\ldots M_{k}\overset{\widehat{\mu }^{\ast }a.s.}{\rightarrow }0$%
, $k\rightarrow \infty $, by Lemma 5.6 of Goldie and Maller (2000) we have $%
I_{-\log
^{-}M_{1}}(\overline{\mathcal{L}(\sum_{i=1}^{L}X_{i})})<\infty $.
The chain of equalities
\begin{equation*}
\infty >I_{-\ln ^{-}M_{1}}(\overline{\mathcal{L}(\sum_{i=1}^{L}X_{i})}%
)=\int_{(1,\infty )}\dfrac{\log x}{\int_{0}^{\ln
x}\widehat{\mu}^{\ast }\{-\log ^{-}M_{1}\leq
-y\}dy}d\widehat{\mu}^{\ast }\{N_{1}\leq x\}=
\end{equation*}
\begin{equation*}
=\int_{(1,\infty )}\dfrac{\log x}{\int_{0}^{\log
x}\widehat{\mu}^{\ast }\{\log M_{1}\leq
-y\}dy}d\widehat{\mu}^{\ast }\{N_{1}\leq
x\}=I_{R_{1}}(\overline{\mathcal{L}(\sum_{i=1}^{L}X_{i})}).\newline
\end{equation*}
finishes the proof. \newline \emph{Case $\alpha\in (0,1)$}. Assume
that Condition D$_{\beta_{1}}$ as well as (d) and (e) of Theorem 2
hold. Put $\alpha:=\beta_{1}$. According to what we have already
proved, the modified transform $\mathbb{T}_{\alpha}$ (defined in
the Introduction) has a fixed point $\mu_{1}$ with mean $m$, say.
Its LST $\varphi_{1}$ satisfies the equality
\begin{equation*}
\varphi_{1} (s)=\mathbb{E}\prod_{i=1}^{L }\varphi_{1}
(X_{i}^{\alpha}s).
\end{equation*}
Set $\varphi_{\alpha}(s):=\varphi_{1}(s^{\alpha})$. The so defined
function is the LST of a distribution $\mu_{\alpha}$ given by (7).
Moreover, $\varphi_{\alpha}(s)$ satisfies (3) and (6). Note that
it is easy to check that in that case the right hand side of (3)
is well-defined. Thus $\mu_{\alpha}$ is the $\alpha$-elementary
fixed point.\newline In the reverse direction, let $\mu_{\alpha}$
be an $\alpha$-elementary fixed point with the LST
$\varphi_{\alpha}$. By Proposition 1, Condition D$_{\beta_{1}}$
holds with $\beta_{1}=\alpha$. It remains to show that the
conditions (d) and (e) are necessary. Define the nonnegative
nondecreasing and continuous function
$\psi_{\alpha}(s):=\varphi_{\alpha}(s^{1/\alpha})$ (it is
reasonable to call the move from $\varphi_{\alpha}$ to
$\psi_{\alpha}$ as \emph{the inverse stable transformation}). It
satisfies the equality
\begin{equation*}
\psi_{\alpha} (s)=\mathbb{E}\prod_{i=1}^{\infty }\psi_{\alpha}
(X_{i}^{\alpha}s),
\end{equation*}
which is the analogue of (3) for the modified transform
$\mathbb{T}_{\alpha}$, and
\begin{equation}
\underset{s\rightarrow
+0}{\lim}\dfrac{1-\psi_{\alpha}(s)}{s}=m\text{.}
\end{equation}
As $\mathbb{T}_{\alpha}$ verifies Condition D$_{1}$, by Lemma 12
$\psi_{\alpha}$ is the LST of a fixed point of
$\mathbb{T}_{\alpha}$. This fixed point has finite mean in view of
(21). Thus, according to the first part of the proof, the
conditions (d) and (e) are indeed necessary. The proof is
complete.\newline \textbf{Proof of Proposition 3.} (a-b) Keeping
in mind the proof of Theorem 2, it remains to show that given
$m>0$ in (6), there exists a unique $\alpha$-regular fixed point.
The proof below is standard, but it is included here for
completeness. Suppose that there exist two $\alpha$-elementary
fixed points whose LST $\varphi_{1,\alpha}$ and
$\varphi_{2,\alpha}$(say) satisfy (6) with the same $m$. From (3),
we deduce that the function $\Xi_{\alpha} (s):=\dfrac{\left|
\varphi _{1,\alpha}(s)-\varphi _{2,\alpha}(s)\right|
}{s^{\alpha}}$
satisfies the inequality $\Xi (s)\leq \mathbb{E}\Xi (B^{(\alpha)}s)$, $%
s> 0$. Iterating this $n$ times gives $\Xi (s)\leq \mathbb{E}\Xi
(\exp(S_{n}^{(\alpha)})s)$, $s>0$. By Theorem 2,
$\exp(S_{n}^{(\alpha)})$ a.s. goes to zero, when
$n\rightarrow\infty$. By Lemma 13, $\varphi_{1,\alpha}\equiv
\varphi_{2,\alpha}$. \newline (c) Assume that there exists a
nonelementary fixed point with the LST $\varphi^{*}$. By
Proposition 1, $1-\varphi^{*}$ is regularly varying at zero with
the index $\alpha$. In view of the assumption, the corresponding
slowly varying function $l$, say, is not equivalent to a constant.
In fact, we have $l(s)\rightarrow+\infty$ when $s\rightarrow +0$.
If $\alpha=1$, it is obvious, if $\alpha\in (0,1)$, use the
inverse stable transformation (defined in the proof of Theorem 2)
and Lemma 12 to reduce this case to the previous one. Now we use
the idea of the proof of Theorem 7.4 of Liu (1998). According to
parts (a) and (b) of the Proposition, for each $m>0$ in (6) there
exists an $\alpha$-elementary fixed point with the LST
$\varphi_{\alpha}^{(m)}$, say. Moreover, there exists an $s_{m}>0$
such that  $1-\varphi_{\alpha}^{(m)}(s)\leq 1-\varphi^{*}(s)$, for
all $0<s<s_{m}$. Now Lemma 7.3 of Liu (1998) applies. This yields
for each $m>0$, $1-\varphi_{\alpha}^{(m)}(s)\leq
1-\varphi^{*}(s)$, for all $s>0$. As
$\varphi_{\alpha}^{(m)}(s)\rightarrow0$ as $m\rightarrow +\infty$,
we get $\varphi^{*}(s)=0$, $s>0$. A contradiction. The proof is
complete.\newline \textbf{Proof of Proposition 4.} We will use the
notation exploited in the
proof of Theorem 2. Set also $\widehat{W}:=\underset{n\rightarrow \infty }{\lim \inf }%
W^{(n)}(1)$. Assume that $\mathbb{E}W^{p}=\mathbb{E}_{\mu }(%
\widehat{W})^{p}<\infty $. Then $\underset{n\rightarrow \infty }{\lim }$ $%
\mathbb{E}_{\widehat{\mu }}\left( W^{(n)}(1)\right) ^{p-1}=\mathbb{E}_{%
\widehat{\mu }}(\widehat{W})^{p-1}<\infty $, the inequality being
implied by Lemma 11. An appeal to (20) reveals that the condition
(8) is necessary. Assume now that (8) holds or which is equivalent
$\mathbb{E}M_{1}^{p-1}<1$
and $\mathbb{E}N_{1}^{p-1}<\infty $. These inequalities ensure that $V_{1}:=%
\underset{n\rightarrow \infty }{\lim \inf }V_{1,n}<\infty $ (by
the general
theory of perpetuities, see Goldie and Maller (2000)) and $\mathbb{E}_{%
\widehat{\mu }^{\ast }}V_{1}^{p-1}<\infty $ (by the triangle inequality in $%
L_{p-1}$).  By Fatou's lemma, we have
\begin{equation*}
\mathbb{E}_{\widehat{\mu }^{\ast }}(V_{1})^{p-1}\geq \mathbb{E}_{\widehat{%
\mu }^{\ast }}\left( \mathbb{E}_{\widehat{\mu }^{\ast }}(\widehat{W}/%
\mathcal{G})\right) ^{p-1}.
\end{equation*}
As it was announced we only consider the case $p\in (1,2]$. In
view of
Jensen's inequality, the right hand side is bounded below by $\mathbb{E}_{%
\widehat{\mu }^{\ast }}\left( \mathbb{E}_{\widehat{\mu }^{\ast }}(\widehat{W}%
^{p-1}/\mathcal{G})\right) =\mathbb{E}_{\widehat{\mu }^{\ast }}\widehat{W}%
^{p-1}$. It remains to use Lemma 11. The proof is
complete.\newline
\textbf{Proof of Proposition 6.} For $\nu _{1}$ and $\nu _{2}\in \mathcal{P}%
^{+}(\delta ,m)$ with characteristic functions $\psi _{1}(s)$ and
$\psi _{2}(s)$ respectively,
let us denote by $\varphi _{i}(s)$ the characteristic functions of $\mathbb{T}\nu _{i}$, $i=1,2$%
. For any complex $z_{i}$, $Z_{i\text{ }}$with $\left|
z_{i}\right| \leq 1$ and $\left| Z_{i}\right| \leq 1$ we have
\begin{equation}
\left| \prod z_{i}-\prod Z_{i}\right| \leq \sum \left|
z_{i}-Z_{i}\right| \text{,}
\end{equation}
when the right hand side is finite. Thus we obtain
\begin{eqnarray*}
\left| \varphi _{1}(s)-\varphi _{2}(s)\right|  &=&\left| \mathbb{E}%
\prod_{i=1}^{\infty }\psi
_{1}(X_{i}s)-\mathbb{E}\prod_{i=1}^{\infty }\psi
_{2}(X_{i}s)\right| \newline
\leq  \\
&\leq &\mathbb{E}\sum_{i=1}^{\infty }\left| \psi _{1}(X_{i}s)-\psi
_{2}(X_{i}s)\right| .
\end{eqnarray*}
Recalling that the rv $B^{(1)}$ is defined in Proposition 1 and setting $%
f(s):=\dfrac{\left| \psi _{1}(s)-\psi _{2}(s)\right| }{s}$, we get
\begin{equation*}
r_{p}(\mathbb{T}\nu _{1},\mathbb{T}\nu _{2})=\int_{0}^{\infty
}s^{-p-1}\left| \varphi _{1}(s)-\varphi _{2}(s)\right| ds\newline
\leq
\end{equation*}
\begin{equation*}
\leq \int_{0}^{\infty }s^{-p}\mathbb{E}\sum_{i=1}^{\infty
}X_{i}f(X_{i}s)ds=\int_{0}^{\infty
}s^{-p}\mathbb{E}f(B^{(1)}s)\newline ds\leq
\end{equation*}
\begin{equation*}
\leq \mathbb{E}(B^{(1)})^{p}\int_{0}^{\infty }z^{-p-1}\left| \psi
_{1}(z)-\psi _{2}(z)\right| dz=
\end{equation*}
\begin{equation*}
=\mathbb{E}\left( \sum_{i=1}^{L}X_{i}\right) ^{p}r_{p}(\nu
_{1},\nu _{2})
\end{equation*}
Among others this justifies using (22). The proof is
complete.\newline
\textbf{Proof of Proposition 7.} Define the function $t(y):=\mathbb{E}%
(\sum_{i=1}^{L}X_{i}^{y})$. It is convex where it is finite. By
convexity, the condition $t(1)=t(b)=1$ implies $t(y)<1$ for $y\in
(1,b )$.
Using Jensen's inequality gives $\mathbb{E}$log$B^{(1)}=\mathbb{E}%
(\sum_{i=1}^{L}X_{i}\log X_{i})<0$. On the other hand, $\mathbb{E}%
(\sum_{i=1}^{L}X_{i})^{b }<\infty $ implies\newline
$\mathbb{E}(\sum_{i=1}^{L}X_{i})\log
^{+}(\sum_{i=1}^{L}X_{i})<\infty $. Therefore, by Theorem 2 with
$\alpha=1$ there exists a fixed point $\mu=\mathcal{L}(W)$ having
finite mean $m>0$, say. \newline The proof goes almost the same
path as that of Proposition 1.2 in Iksanov and Kim (2003), where
the tail bahaviour of fixed points of the (shifted Poisson) shot
noise transforms has been studied. Keeping this in mind, we only
give a sketch of the proof and refer the interested reader to
Iksanov and Kim (2003) for details. There is a point of essential
difference between this work and the just cited one. In place of
the simple perpetuity (24) investigated in Iksanov and Kim (2003)
we should use the other perpetuity (23) found by Liu (2000, Lemma
4.1). At this point we would like to stress that these
perpetuities are quite different and (23) is not reduced to (24)
even for the shot noise transforms.

We start with the renewal equation
\begin{equation*}
P_{b }(x)=\int_{-\infty }^{\infty }P_{b }(x-y)\rho _{b }(dy)+Q_{b
}(x)\text{,}
\end{equation*}
where
\begin{equation*}
P_{b }(x):=e^{-x}\int_{0}^{\exp x}y^{b }\mu (y,\infty )dy\text{, }
\end{equation*}
\begin{equation*}
Q_{b }(x):=e^{-x}\int_{0}^{\exp x}y^{b }(\mu (y,\infty
)-N(y,\infty ))dy
\end{equation*}
and $\rho _{b }(dy):=e^{(b-1)y}\mathcal{L}(B^{(1)})(de^{y})$. We
have the equality of distributions which essentially shows that
$\overline{W}$ with $\mathcal{L}(\overline{W})=\overline{\mu}$ is
the perpetuity
\begin{equation}
\overline{W}\overset{d}{=}B^{(1)}\overline{W}+C\text{,}
\end{equation}
where $\overline{W}$ is independent of ($B^{(1)},C$) and $\mathbb{E}f(C)=%
\mathbb{E}(\sum_{i=1}^{L}X_{i}f(\sum_{\substack{ 1\leq k\leq L  \\ k\neq i}}%
X_{k}))$, for any nonnegative Borel functions $f$.\newline
Set $I(x):=\int_{0}^{\infty }y^{x}(\mu (y,\infty )-N(y,\infty ))dy$. For $%
\beta <b-1$, we have $\mathbb{E(}B^{(1)})^{\beta }=\mathbb{E}%
(\sum_{i=1}^{L}X_{i}^{\beta +1})<1$ and, by Proposition 4 $\mathbb{E}%
W^{\beta +1}<\infty $ which imply
\begin{equation*}
0<I(\beta )=(\beta +1)^{-1}\mathbb{E}W^{\beta
+1}(1-\mathbb{E(}B^{(1)})^{\beta })<\infty .
\end{equation*}
Applying the $c_{\beta }$-inequality to (23) results in
\begin{equation*}
m^{-1}\mathbb{E}W^{\beta +1}(1-\mathbb{E(}B^{(1)})^{\beta })\leq
(2^{\beta -1}\vee 1)M_{\beta },
\end{equation*}
\begin{equation*}
m^{-1}\mathbb{E}W^{\beta +1}(1-\mathbb{E(}B^{(1)})^{\beta })\geq
(2^{\beta -1}\wedge 1)K_{\beta },
\end{equation*}
where the constant
\begin{equation*}
M_{\beta }:=\mathbb{E}(\sum_{i=1}^{L}X_{i})^{\beta +1} \text{ , if
\ } \beta \in (0,1] \text { \ \ and \ } :=\mathbb{E}W^{\beta
}\mathbb{E}(\sum_{i=1}^{L}X_{i})^{\beta +1} \text{, otherwise,}
\end{equation*}
is identified in Lemma 4.2 of Liu (2000), and the
constant
\begin{equation*}
K_{\beta }:=\mathbb{E}(\sum_{i=1}^{L}X_{i})^{\beta +1} \text{ , if
\ } \beta \geq 1 \text { \ \ and \ } :=\mathbb{E}W^{\beta
}\mathbb{E}(\sum_{i=1}^{L}X_{i})^{\beta +1} \text{, otherwise,}
\end{equation*}
can be obtained in the similar way. Letting $\beta $ go to $b-1$
along some subsequence gives
\begin{equation*}
0<mb^{-1}(2^{b-2}\wedge 1)K_{b-1}\leq I(b-1)\leq mb
^{-1}(2^{b-2}\vee 1)M_{b-1}<\infty \text{,}
\end{equation*}
which according to Lemma 9.2 of Goldie (1991) implies that
$Q_{b}(x)$ is directly Riemann integrable.

By the key renewal theorem for the whole line we obtain: \newline
1) if $\mathcal{L}(\log B^{(1)})$ is nonarithmetic then
\begin{equation*}
\underset{x\rightarrow \infty }{\lim }P_{b}(x)=(\mathbb{E}%
(\sum_{i=1}^{L}X_{i}^{b}\log X_{i}))^{-1}\int_{0}^{\infty
}y^{b-1}(\mu (y,\infty )-N(y,\infty ))dy:=C_{b };
\end{equation*}
\newline
2) if $\mathcal{L}(\log B^{(1)})$ is arithmetic with the span
$\varsigma $ then for all $x\in \mathbb{R}$\newline
\begin{equation*}
\underset{n\rightarrow \infty }{\lim }P_{b }(x+\varsigma n)=(\mathbb{E}%
(\sum_{i=1}^{L}X_{i}^{b}\log X_{i}))^{-1}\sum_{k=-\infty }^{\infty
}Q_{b}(x+\varsigma k):=C_{b }(x)\text{. }
\end{equation*}
Thus we only need to consider the arithmetic case. \newline From
the results of Grincevi\v{c}ius (1975) it follows that there exist
$d_{1}>0$ and $d_{2}<\infty $ such that, for all $x$ large enough,
we have
\begin{equation*}
0<d_{1}\leq x^{b }\mu (x,\infty )\leq d_{2}<\infty \text{.}
\end{equation*}
Taking this into account, we can prove that $G_{b }(x)$ slowly
varies at $\infty $. The working from Iksanov and Kim (2003) may
be repeated, but we must have an alternative proof of the fact
that $\underset{n\rightarrow \infty }{\lim }\dfrac{N(t_{n}u,\infty
)}{\mu (t_{n}u,\infty )}\in \lbrack
0,1]$ provided the limit exists for some sequence $t_{n}\rightarrow \infty $%
, as $n\rightarrow \infty $. This follows easily from (23) which
implies that $\overline{\mu }(x,\infty )\geq N^{\ast }(x,\infty
)$, thus giving $\mu (x,\infty )\geq N(x,\infty )$ (use the LST).
It remains to show that
\begin{equation*}
\underset{n\rightarrow \infty }{\lim }e^{(x+\varsigma n)}\mu
(e^{x+\varsigma
n},\infty )=C_{b }(0)\text{ locally uniformly in }e^{x}\text{ on }%
(0,\infty )\text{,}
\end{equation*}
which together with slow variation will give the result.

\section{A Pitman-Yor problem as a particular case}

In Pitman and Yor (2000, p.35) the following problem concerning
distributions on the nonnegative half line has been mentioned. For
what distributions $\nu $ there exists a distribution $\mu $ with
finite mean such that
\begin{equation}
\overline{Y}\overset{d}{=}A\overline{Y}+Y\text{,}
\end{equation}
where\ $\mathcal{L}(A)=\nu $, $\mathcal{L}(Y)=\mu $ and $\mathcal{L}(%
\overline{Y})=\overline{\mu }$.

To point out the solution to the above problem, let us introduce
the random walk $T_{0}:=0$, $T_{n}:=\sum_{k=1}^{n}\log A_{k}$,
$k=1,2,\ldots$, where $A_{1}$, $A_{2},\ldots $ are independent
copies of the rv $A$. We now intend to explain how the next result
originally obtained in Iksanov and Kim (2003) follows from Theorem
2 (case $\alpha=1$).\newline \textbf{Proposition 8. }The condition
\begin{equation}
\underset{n\rightarrow \infty}{\lim} T_{n}=-\infty \text{ \ almost
surely}
\end{equation}
is necessary and sufficient for the existence of a distribution
$\mu \neq \delta _{0}$ satisfying (24). Given $m>0$, there exists
a unique $\mu $ with mean $m$.

Recall that a right-continuous and nonincreasing function
$g:(0,\infty )\rightarrow \lbrack 0,\infty )$ allows us to define
the generalized inverse function $g^{\leftarrow }$ as follows
\begin{equation*}
g^{\leftarrow }(z):=\inf \{u:g(u)<z\}\text{, if }z<g(0^{+})\text{ ; }:=0\text{%
, otherwise.}
\end{equation*}

We do not know how the problem could be solved on using only the
equality (24). Hence, the alternative representation of the rv $Y$
is given next.
\newline \textbf{Lemma 9. }(Iksanov and Kim (2003)\textbf{)} a)
Assume that a distribution of the rv $Y$ with finite mean $m>0$
satisfies (24), where the rv $A$ is such that\newline $\gamma
:=\mathbb{P}\{A=0\}\in \lbrack 0,1).$ Then we have
\begin{equation}
Y\overset{d}{=}m\gamma +\sum_{i=1}^{\infty }Y_{i}h(\tau
_{i})\text{, }
\end{equation}
where $Y_{1},Y_{2},\ldots$ are independent copies of $Y$, which
are also
independent of a Poisson flow $\{\tau _{i}\},i\geq 1$ with intensity $1$; $%
h:(0,\infty )\rightarrow \lbrack 0,\infty )$ is right-continuous
and nonincreasing function defined by
\begin{equation*}
h^{\leftarrow }(x)=\int_{x}^{\infty
}z^{-1}\mathcal{L}(A)(dz),x>0\newline \text{.}
\end{equation*}
Therefore, we have $\int_{0}^{\infty }h(z)dz=1-\gamma $.\newline
b) If a distribution of $Y$ satisfies (26), where the intensity of
the Poisson flow is equal to $\lambda >0$, and a right-continuous
and nonincreasing function $h:(0,\infty )\rightarrow \lbrack
0,\infty )$ satisfies the equality
\begin{equation*}
\lambda \int_{0}^{\infty }h(z)dz=1-\gamma \text{,}
\end{equation*}
Then $\mathcal{L}(Y)$ solves (24), where $\mathcal{L}(A)$ is defined by $%
\mathcal{L}(A)(dx)=-\lambda xh^{\leftarrow }(dx)$.\strut

We are ready to check that the assertion of Proposition 8
with\newline $\mathbb{P}\{A=0\}=0$ is contained in Theorem 2. To
this end, we use the observation that $\mathcal{L}(Y)$ solving
(24) is a fixed point of the Poisson shot noise transform (see
Iksanov and Yurek (2002) for more details) or, in other words,
$\mathcal{L}(Y)$ satisfies (26) and vice versa. Hence, let us set
in Theorem 2 $X_{i}=h(\tau _{i})$, where $h$ and $\tau _{i}$ are
defined in Lemma 9.

Condition D$_{1}$ is equivalent to $\int_{0}^{\infty} h(u)du=1$
(when applying part b) of Lemma 9 we can always take a Poisson
flow with unit intensity). Thus, all that we need is to show that
the condition (d) of Theorem 2 implies the condition (e) of the
same Theorem. Under the current notation, this reduces to showing
that (25) implies
\begin{equation}
I_{\log A}(\overline {\mathcal{L}(\sum_{i=1}^{\infty}h(\tau
_{i})))}<\infty.
\end{equation}
To see this, it suffices to note that
$\chi(dt)=-h^{\leftarrow}(dt)$ and
$\mathcal{L}(A)=\chi^{\ast}$.\newline Condider two cases: (1)
$-\infty<\mathbb{E}\log A<0$ and (2) $\mathbb{E}\log A=-\infty$ or
$\mathbb{E}\log A$ does not exist.\newline \emph{Case 1.} Our task
simplifies to checking that
\begin{equation}
\mathbb{E}\log A\in(-\infty, 0) \text{ \ implies \ }
\mathbb{E}\log(1+\overline
{\sum_{i=1}^{\infty}h(\tau_{i}))}<\infty.
\end{equation}
We have $\mathbb{E}\log A\in(-\infty, 0)$ implies
$\mathbb{E}\log(1+ A)<\infty$ and Condition D$_{1}$
($\mathbb{E}\sum_{i=1}^{\infty}h(\tau_{i})=1$) implies
$\mathbb{E}\log(1+\sum_{i=1}^{\infty}h(\tau_{i}))<\infty$. It is
easy to observe that $\overline
{\sum_{i=1}^{\infty}h(\tau_{i}))}\overset{d}{=}\sum_{i=1}^{\infty}h(\tau_{i})+A$
(use LST's). Therefore, applying the inequality (29) gives
(28).\newline\emph{Case 2.} Set
\begin{equation*}
g(x):=\dfrac{x}{\int_{0}^{x}\mathbb{P}\{\log A\leq -y\}dy},\
f(x):=xg(\log (1+x)), x\geq 0.
\end{equation*}
\begin{equation*}
t(x):=(x\wedge 1)(1+g(\ln (1+x))), x\geq 0.
\end{equation*}
Since $g(x)/x$ is nonincreasing, we have $g(x+y)\leq g(x)+g(y)$,
for all $x,y\geq 0$ (subadditivity). This together with
nondecreasingness of $g(x)$ and the inequality
\begin{equation}
\log (1+x+y)\leq \log (1+x)+\log (1+y),
\end{equation}
results in
\begin{equation*}
f(x+y)/(x+y)\leq f(x)/x+f(y)/y \text{ \ , for all \ } x,y\geq 0.
\end{equation*}
Consequently, the function $t(x)$ is submultiplicative. By Theorem
25.3 of Sato (1999) the integrability of a submultiplicative
function with respect to an infinitely divisible distribution is
equivalent to the integrability (near infinity) of the function in
question with respect to the corresponding L\'{e}vy measure. As
$(-1)h^{\leftarrow }(dx)$ defines the L\'{e}vy measure of the
infinitely divisible distribution $\mathcal{L}(\sum_{i=1}^{\infty
}h(\tau _{i}))$, we have
\begin{equation*}
\int_{0}^{\infty }t(x)\mathcal{L}(\sum_{i=1}^{\infty }h(\tau
_{i}))(dx)<\infty \text { \ if and only if} \int_{1}^{\infty
}t(x)(-1)h^{\leftarrow }(dx)<\infty.
\end{equation*}
By the same criterion, Condition D$_1$ implies that both integrals
$\int_{1}^{\infty }xh^{\leftarrow }(dx)$ and $\int_{0}^{\infty }x%
\mathcal{L}(\sum_{i=1}^{\infty }h(\tau _{i}))(dx)$
 converge. This and the above display together implies
\begin{equation*}
\infty>\int_{1}^{\infty }f(x)\mathcal{L}(\sum_{i=1}^{\infty
}h(\tau _{i}))(dx) \text { \ if and only if \ } \infty>
\int_{1}^{\infty }f(x)(-1 )h^{\leftarrow }(dx)=.
\end{equation*}
\begin{equation*}
=\int_{1}^{\infty }g(\log(1+x))d\mathbb{P}\{A\leq x\}\text { \ if
and only if \ } \infty>\int_{0}^{\infty
}g(\log^{+}x)d\mathbb{P}\{A\leq x\}
\end{equation*}
It is clear that the former inequality is equivalent to (27).
Recall that we consider the case when $\log A$ has no finite mean.
Thus, by (1.19) of Kesten and Maller (1996), the latter inequality
is equivalent to (25), and the asserted follows.

\section{Comments and some references}5.1. Fixed
points of the smoothing transforms with a.s. finite number of
summands have been receiving much attention in the literature. It
is worth mentioning that usually more or less restrictive
additional moment assumptions have been imposed. The basic
techniques and results for the case of finitely many summands were
developed in the basic paper of Durrett and Liggett (1983). Liu
(1998) extended their results to the case of finite but random
number of summands. Fixed points of the BRW smoothing transforms
with almost surely finite number of summands were studied (without
using the name) by Biggins (1977), Biggins and Kyprianou (1996,
1997, 2001, 2003), Kyprianou (1998), Liu (2000). Except for the
papers on fixed points of the shot noise transforms (see point 5.2
below), we are aware of only four works dealing with fixed points
of the smoothing transforms with \emph{infinite} number of
summands. These are Lyons (1997), R\"{o}sler (1992), Caliebe and
R\"{o}sler (2003a,b). Note also that the last three papers
investigate fixed points concentrated on the whole line. Since the
extensive surveys were presented, for example, in R\"{o}sler
(1992) and Liu (1998), we refrain from listing other works on the
subject here.
\newline 5.2. Details regarding \emph{some} fixed points of Poisson shot noise
transforms can be found in Iksanov (2002a,b), Iksanov and Jurek
(2002), Iksanov and Kim (2003, 2004). The Pitman-Yor problem was
solved in Iksanov and Kim (2003) by using an approach different
from that taken here.
\newline 5.3. In the area of smoothing transforms the notion
of "stable transformation" is due to Durrett and Liggett (1983).
See also Guivarc'h (1990) and Liu (1998) for further development
of this concept.
\newline 5.4. The trick based on the martingale stopping we have
used in the proof of Proposition 1 had come to our attention from
Durrett and Liggett (1983, Theorem 2.18). Similar idea, but
sometimes in different forms, was much exploited in works
connected with convergence of Markov processes and/or martingales.
Here we only mention Biggins and Kyprianou (1997, Sections 6 and
8; 2001; 2003, Lemma 1) who developed approaches of such a flavour
in relations to fixed points of the smoothing transforms.

\section{Appendix} For ease of references we collect here
some facts taken mainly from other sources. \newline \textbf{Lemma
10.} (a) (Liu (1998), Lemma 1.1) If $\mathbb{P}\{X_{i}=0$ or $1$,
for all $i\leq L\}=1$ then there are no fixed points.\newline (b)
(Liu (1998), Lemma 3.1) If fixed points exist then
$\mathbb{E}L>1$.\newline \textbf{Lemma 11.} Let $\theta $ be a
finite measure and $\nu $ a
probability measure on a $\sigma $-field $\mathcal{F}$. Suppose that $%
\mathcal{F}_{n}$ are increasing sub-$\sigma $-fields whose union generates $%
\mathcal{F}$ and that the restriction of $\theta $ to
$\mathcal{F}_{n}$ is
absolutely continuous with respect to the restriction of $\nu $ to $\mathcal{%
F}_{n}$ with Radon-Nikodym derivative $Y_{n}$. If
$Y:=\underset{n\rightarrow \infty }{\lim \sup }Y_{n}<\infty $ then
for fixed $\varepsilon \geq 0$ $\int X^{1+\varepsilon }d\nu
<\infty $ if and only if $\int X^{\varepsilon }d\theta <\infty $.
If one of these is finite then $\int X^{1+\varepsilon }d\nu =$
$\int X^{\varepsilon }d\theta $.\newline \textbf{Proof.
}When\textbf{\ }$\varepsilon =0$ this is just Lemma 10.2 from
Lyons and Peres (2004) where among others it was shown that
$\theta (A)=\int Xd\nu $ for all $A\in \mathcal{F}$. Hence,
$X^{1+\varepsilon }d\nu =X^{\varepsilon }d\theta $ and the Lemma
follows.

The main message of the next Lemma is that the stable and inverse
stable transformations give a one-to-one correspondence between
fixed points of $\mathbb{T}$ with
$\mathbb{E}\sum_{i=1}^{L}X_{i}^{\beta_{1}}=1$, $\beta_{1}\in
(0,1)$ and those of the modified transform
$\mathbb{T}_{\beta_{1}}$. For the elementary fixed points this
observation was exploited in the proof of the $\alpha\in (0,1)$
case of Theorem 2. \newline \textbf{Lemma 12.} Assume that
Condition D$_{\beta_{1}}$ holds and a nonnegative nonincreasing
and continuous function $f(s), s\geq 0$ satisfies $f(0)=1$ and
\begin{equation*}
f(s)=\mathbb{E}\prod_{i=1}^{L}f(X_{i}s)\text{,}
\end{equation*}
and $1-f(s)$ regularly varies at zero with index
$\beta_{1}$.\newline Then $f$ is the LST of a fixed point of
$\mathbb{T}$.\newline \textbf{Proof.} The proof follows the
similar path as that of part (d) of Theorem 7.1 of Liu (1998). A
close inspection of Liu's proof reveals that his assumption that
$\phi$ is the LST of a fixed point is not needed. (Also the
restrictions $L<\infty$ and (H1) are of no importance for the
stated \emph{here} result to hold). It simply suffices to require
that $\phi$ satisfies the conditions of our Lemma.\newline Indeed,
assume first that $\beta_{1}\in (0,1)$. First we want to show how
one may construct an LST $g$, say, such that
\begin{equation*}
\underset{s\rightarrow +0}{\lim}\dfrac{1-f(s)}{1-g(s)}=1\text{.}
\end{equation*}
Even though $1-f$ is regularly varying, it is far from being
obvious to us that such functions $g$ do exist. By Theorem 1.7.6
and a variant of Theorem 1.5.8 of Bingham, Goldie and Teugels
(1989),
\begin{equation*}
h(s):=(\alpha/\Gamma(1-\alpha))\int_{0}^{s} \int_{0}^{\infty}
e^{-ux}(1-f(1/x))dxdu
\end{equation*}
is the nonnegative function with completely monotone derivative
and\newline $\underset{s\rightarrow
+0}{\lim}\dfrac{1-f(s)}{h(s)}=1$. If $\underset{s\rightarrow
\infty}{\lim} h(s)=A\in (1,+\infty]$ then there exists a finite
nonzero $s_{1}$ such that $h(s_{1})=1$. In that case, the
function\newline $g(s):=1-h(s_{1}(1-e^{-s/s_{1}}))$ is a wanted
LST. Indeed, the superposition of two nonnegative functions with
completely monotone derivatives is a nonnegative function with
completely monotone derivative and $1-e^{-s}\sim s$ as
$s\rightarrow 0$. If $\underset{s\rightarrow \infty}{\lim}
h(s)=A\in (0,1]$ then choose $g(s):=1-h(s)$.

Lemma 7.2 of Liu (1998) applies since his condition
$\rho(\alpha)\leq 1$ is implied by our Condition  D$_{\beta_{1}}$.
Hence, the statement of Lemma 7.3 of the same reference is true in
our case too. Thus we have that $f$ is completely monotone as the
pointwise limit of the sequence of the LSTs
\begin{equation*}
\psi _{0}(s):=g(s)\text{, \ \ }\psi
_{n}(s):=\mathbb{E}\prod_{i=1}^{L}\psi_{n-1}(sX_{i})\text{,}\\
n=1,2,\ldots
\end{equation*}
Since by assumption $f(0)=1$, $f(s)$ is the LST of a
(nondefective) distribution and therefore it is the LST of a fixed
point. Assume now that $\beta=1$. Unfortunately, it is not clear
to us how we could construct a function like $h$ in this case.
Thus we will use another argument. For each $\gamma\in (0,1)$,
consider the collection of point processes $Z_{\gamma}(\cdot)$
with the points $\{X_{\gamma,i}\}_{i=1}^{L}$ such that
$X_{\gamma,i}:=X_{i}^{1/\gamma}$, $\overline{i=1,L}$. Define also
$f_{\gamma}(s):=f(s^{\gamma})$. By what we have already proved,
for each $\gamma\in (0,1)$, $f_{\gamma}(s)$ are the LST's of the
BRW smoothing transform constructed by the point process
$Z_{\gamma}(\cdot)$. As $f(s)=\underset{\gamma\rightarrow 1}{\lim}
f_{\gamma}(s)$ and $f(0)=1$, $f(s)$ is the LST of a fixed point.
The proof is complete.

The following Lemma contains an observation implicitly made in
Theorem 1 of Athreya (1969). It provides an easy way of showing
that when restricted to the set of the $\alpha$-elementary fixed
points, fixed points are unique up to the scale (see the proof of
Proposition 3).\newline \textbf{Lemma 13.} Let $\varphi
_{1,\alpha}(s)$ and $\varphi _{2,\alpha}(s)$ be the
Laplace-Stieltjes transforms of distributions and assume that they
satisfy (6) with the same $m$. If for large enough positive
integer $n$, the function $\Xi_{\alpha} (s):=\dfrac{\left| \varphi
_{1,\alpha}(s)-\varphi _{2,\alpha}(s)\right| }{s^{\alpha}}$
satisfies the inequality
\begin{equation*}
\Xi (s)\leq \mathbb{E}\Xi (C_{n}s)\text{, } s>0,
\end{equation*}
where $\{C_{n}\}$ is a sequence of rvs that goes to zero almost
surely, as $n\rightarrow\infty$, then $\varphi
_{1,\alpha}(s)\equiv $ $\varphi _{2,\alpha}(s)$.\newline
\textbf{Lemma 14.} A fixed point with unit mean exists
if and only if the nonnegative martingale $W^{(n)}(\gamma)$, $%
n=1,2,\ldots $ given by (1) converges in mean to it.\newline
\textbf{Proof.} One implication is obvious, hence let us assume
that a fixed point $W(\gamma)$ with unit mean exists. Let
$\varphi$ be its LST. Arguing as in the proof of Lemma 12, we have
that
$\varphi(s)=\underset{n\rightarrow\infty}{\lim}\varphi_{n}(s)$,
where
\begin{equation*} \varphi _{0}(s):=e^{-s}\text{, \ \ }\varphi
_{n}(s):=\mathbb{E}\prod_{i=1}^{L}\varphi_{n-1}(sX_{i})\text{,}\\
n=1,2,\ldots
\end{equation*}
Since $\varphi _{n}(s)=\mathbb{E}e^{-sW^{(n)}(\gamma)}$,
$W^{(n)}(\gamma)$ weakly converges to $W(\gamma)$ when
$n\rightarrow\infty$. Thus $W^{(n)}(\gamma)$ cannot converge to
zero almost surely. Therefore, it must converge to $W(\gamma)$ in
mean (see the last paragraph on page 3). The proof is complete.

Note that the above result (with different proof) is also given in
Theorem 2.2(1) of Caliebe and R\"{o}sler (2003).
\begin{center}
{\bf {\large Acknowledgement}}
\end{center}
My sincere thanks go to Dr Amke Caliebe and Dr Andreas Kyprianou
whose detailed comments on the manuscript have led to the
essential improvement of presentation.

\end{document}